\documentclass{article}

\usepackage[final,nonatbib]{template}

\usepackage[utf8]{inputenc} % allow utf-8 input
\usepackage[T1]{fontenc}    % use 8-bit T1 fonts
\usepackage{hyperref}       % hyperlinks
\usepackage{url}            % simple URL typesetting
\usepackage{booktabs}       % professional-quality tables
\usepackage{amsfonts,amsmath,dsfont}       % blackboard math symbols
\usepackage{nicefrac}       % compact symbols for 1/2, etc.
\usepackage{microtype}      % microtypography
\usepackage{xcolor}         % colors

\newtheorem{theorem}{Theorem}

\newtheorem{lemma}{Lemma}

\title{\vspace{-0.1in}Catoni-style Confidence Sequences under Infinite Variance\vspace{-0.1in}}

\author{%
 Sujay Bhatt, Guanhua Fang, Ping Li \\
 Cognitive Computing Lab, Baidu Research\\
10900 NE 8th St. Bellevue, WA 98004
  \And
Gennady Samorodnitsky \\
 School of ORIE, Cornell University\\
220 Frank T Rhodes Hall, Ithaca,  NY 14853
}

\newcommand{\bbr}{\mathbb{R}}

\newcommand{\bma}{\begin{matrix*}[r]}
\newcommand{\ema}{\end{matrix*}}

\newcommand{\vep}{\varepsilon}

\begin{document}

\maketitle

\begin{abstract}\vspace{-0.1in}
In this paper, we provide an extension of confidence sequences for settings where the variance of the data-generating distribution does not exist or is infinite. Confidence sequences furnish confidence intervals that are valid at arbitrary data-dependent stopping times, naturally having a wide range of applications. We first establish a lower bound for the width of the Catoni-style confidence sequences for the finite variance case to highlight the looseness of the existing results. Next, we derive tight Catoni-style confidence sequences for data distributions having a relaxed bounded~$p^{th}-$moment, where~$p \in (1,2]$, and strengthen the results for the finite variance case of~$p =2$. The derived results are shown to better than confidence sequences obtained using Dubins-Savage inequality. \vspace{-0.1in}
\end{abstract}

\vspace{-0.1in}
\section{Introduction}
\vspace{-0.05in}

Sequential design of experiments is a classical framework in statistical sampling theory, in which the size and the composition of samples are not fixed in advance and are allowed to be functions of the observations themselves~\cite{Rob52}. Confidence sequence is one  particular tool in sequential design that facilitates anytime-valid inference~\cite{DH67,JJ18,HRM21}. In particular, confidence sequence is a sequence of confidence intervals that is valid at data-dependent stopping times.

Formally, let~$X_1,X_2,\cdots$ be an i.i.d stream drawn from distribution~$P$. The basic object of interest is the unknown mean of this distribution, namely,~$\mu = \mathbf{E}_P[X_1]$. A crude way to quantify the uncertainty associated with the mean estimation problem is via \textit{confidence intervals}. Here one constructs a~$\Sigma(X_1,\cdots,X_t)-$measurable interval~$\texttt{CI}_t$ for each~$t \in \mathbb{N}^+$ such that~$\forall~t \in \mathbb{N}^+$, the following holds:~$\mathbb{P}(\mu \in \texttt{CI}_t) \geq 1- \alpha$, for some coverage  probability~$1-\alpha \in (0,1)$. However, as argued in~\cite{WR22}, confidence intervals may undercover at stopping times. Also, it is well known that the data-dependent peeking at confidence intervals inflates the Type-1 error~\cite{JKP17,HRM21}. This motivates confidence sequences, which provide an universal quantification over~$t$. For a confidence parameter~$\alpha \in (0,1)$ and~$t \in \mathbb{N}^+$, the sequence of random intervals~$\{\texttt{CI}_t\}$ that satisfy~$\mathbb{P}(\forall~t \in \mathbb{N}^+,~\mu \in \texttt{CI}_t) \geq 1- \alpha$ are called a~$(1-\alpha)-$\textit{confidence sequence}.

%Of course, one could adopt a näive approach by defining a confidence interval for every~$t$ and taking a union bound. However, to maintain a total probability of error of~$\alpha$, the confidence intervals must have error probabilities whose sum is bounded by~$\alpha$. This results in a confidence interval width that is larger by a logarithmic factor.

\textbf{Literature.}
Confidence sequences (CS) are instrumental in modern application tools like multi-armed bandits~\cite{ZHH21}, A/B testing~\cite{HRM21}, causal inference~\cite{MC21}, etc. Given the importance, it is not surprising that a significant research effort has been allocated to construct confidence sequences under various distributional assumptions on the data. \cite{DH67,JT89} consider~$P$ as normal distributed and construct confidence sequences, while~\cite{JPL15,JKP17} consider~$P$ belonging to an exponential family. \cite{WR20} consider arbitrary but bounded~$P$, while~\cite{HRM21} consider~$P$ having a bounded moment generating function.

Recently,~\cite{WR22} make a substantial contribution to the literature on confidence sequences by relaxing the distributional assumptions to requiring only the existence of a bounded second moment. This is made possible by using robust mean estimator developed in~\cite{Cat12}, which uses influence functions to stabilize the effect of the outliers. While the previous works~\cite{HRM21,WR20} required a Chernoff-type assumption on the distribution resulting in~$O(\sqrt{\log t/ t})$ shrinkage rates for the confidence sequences,~\cite{WR22} show that employing Catoni's estimator improves the rate to~$O(\sqrt{\log \log 2t/ t})$ under weaker assumptions on the distribution. The significance of this result is that there is no excess compromise in considering confidence sequences or weakening the distributional assumptions. A careful reading of~\cite{WR22} shows however that there are inefficiencies compared to the lower bound that can be improved. Also,~\cite{WR22} cannot deal with situations where even the second moment does not exist.

The key improvements over the existing results are as follows:
\begin{itemize}
    \item We extend the analysis to the case of~$p \in (1,2)$, thereby relaxing the distributional assumptions on the data. This allows sequential design for a larger class of distributions, while achieving good control over the width of the CS.
    \item We then derive lower bounds for the width of Catoni-style confidence sequences (CS) for the finite variance case. This shows that~\cite{WR22} almost - \textit{but not quite} - matches the lower bound, leaving room for improvement.
    \item We improve on the näive CS derived using generalized Dubins-Savage inequality in terms of the growth of the width of the confidence interval (CI) w.r.t the confidence parameter~$\alpha$. We improve on the Catoni-style CS in~\cite{WR22} by obtaining a smaller CI width.
\end{itemize}

Note that we consider i.i.d data as it has much wider applicability and has simpler notation, however, the results carryover for stochastic processes with constant conditional expectation using the influence function and the supermartingale arguments as in~\cite{WR22}.

\textbf{Organization.} The paper is organized as follows. Sec.~\ref{sec:PO} formulates the problem and the key objectives that we seek to address.
Sec.~\ref{sec:DS} presents a first attempt at deriving confidence sequences using the well known Dubins-Savage inequality for~$p \in (1,2]$. Sec.~\ref{sec:LB} derives the lower bound for the width of the confidence sequences when the data has finite variance. Sec.~\ref{sec:IV} presents the improvement for the finite variance case and the generalization to the infinite variance case.

\vspace{-0.05in}
\section{Problem Formulation}\label{sec:PO}
\vspace{-0.05in}

In this section, we first introduce the notation and describe the main problem considered in this paper. Recall that~$X_1,X_2,\cdots$ are a sequence of random variables with mean~$\mathbb{E}[X_1] = \mu$ and the $p^{th}-$moment $\mathbb{E}|X_1 - \mu|^p \leq \upsilon_p$ for~$p\in(1,2]$. A level~$1-\alpha \in (0,1)$ confidence sequence for~$\mu$ is a sequence of real numbers~$L_n(X_1,\cdots,X_n)$ and $U_n(X_1,\cdots,X_n)$, where~$L_n,U_n: \mathbb{R}^n \rightarrow \mathbb{R}$ with~$n = 1,2,\cdots$ such that~$L_n \leq U_n$ point-wise and
\[
\mathbb{P}\Big(L_n(X_1,\cdots,X_n) \leq \mu \leq U_n(X_1,\cdots,X_n),  \forall n\geq 1  \Big) \geq 1 - \alpha.
\]

Let the width of the $n^{th}$ confidence interval
\[
W_n(X_1,\cdots,X_n) := U_n(X_1,\cdots,X_n) - L_n(X_1,\cdots,X_n).
\]
\textbf{Objective.} We are interested in how fast this width~$W_n$ can shrink as~$n$ increases, when the $p^{th}-$moment of the data distribution~$P$ is bounded with~$p\in(1,2]$.

\vspace{-0.05in}
\section{Confidence Sequence for Infinite Variance via Dubins-Savage} \label{sec:DS}
\vspace{-0.05in}

In this section, we take the first steps to derive confidence sequences for the infinite variance case. We borrow the width optimization ideas from~\cite{WR20,WR22} and make use of the~$L_p$ version of the classical Dubins-Savage. A generalization of the classical Dubins-Savage inequality~\cite{DS65} first appeared in~\cite{Kal75}, and most recently in~\cite{Kha09}. Lemma~\ref{lem:DS} is a restatement and can be derived using an approach of Doob~\cite{Kha09}.

\begin{lemma} \label{lem:DS}
Let~$\{X_t\}$ be a real-valued stochastic process adapted to the filtration~$\{\mathcal{F}_t\}$, where~$\mathcal{F}_0$ is the trivial sigma-algebra. Let~$\{S_t\}$ be a martingale with~$V_t = \mathbb{E}[|S_t - S_{t-1}|^{p}| \mathcal{F}_{t-1}]$. Then for all~$a\geq 0,~b>0$, we have~$\mathbb{P}\Big(S_t \geq a + b \sum_{i=1}^t V_i \Big) \leq \frac{1}{(1+m_p a b^{\frac{1}{p-1}})^{p-1}},$
where $m_p = \Big(\frac{p-1}{2^{2-p}} \Big)^{\frac{1}{p-1}}$.
\end{lemma}

We have the following result for the confidence sequence.

\begin{theorem}[Dubins-Savage] \label{thm:DS}
 Let~$a = \frac{1}{m_p b^{\frac{1}{p-1}}} \Big(\Big(\frac{2}{\alpha} \Big)^{\frac{1}{p-1}} -1\Big)$. The width of the~confidence interval using Dubins-Savage inequality is given as~
 \[
 W_n(X_1,\cdots,X_n) = \frac{2 b \upsilon_p \sum_{i=1}^n \lambda_i^p}{\sum_{i=1}^n \lambda_i},\]
where the coefficients~$\lambda_i$ are chosen as $\lambda_t = \Big(\frac{a} {t b \upsilon_p (p-1)} \Big)^{1/p}$ for $t \leq n$.
\end{theorem}
From Theorem~\ref{thm:DS}, the width of the confidence interval shrinks at the rate~$O(\frac{\log t}{t^{\frac{p-1}{p}}})$. When~$p=2$, we obtain~$\widetilde{O}(t^{-1/2})$ which is known to be unimprovable~\cite{HRM21}. The dependence on~$\alpha$ is~$O(\alpha^{-1/p})$.
\vspace{-0.05in}
\section{Lower Bound for Finite Variance ($p=2$)} \label{sec:LB}
\vspace{-0.05in}

In this section, we establish the lower bounds for the width of Catoni-style confidence sequences by making similar assumptions as in~\cite{WR22}. We make use of the general law of iterated logarithm by~\cite{Wit85} to support key arguments.

\cite{WR22} discuss one possible way to construct a Catoni-style confidence sequence and is as follows. Suppose the observations~$X_i$ have a finite second moment and~$\texttt{Var}(X_1) \leq \sigma^2$ for a known~$\sigma^2$. Consider an increasing function~$\psi$ such that for~$x \in \mathbb{R}$ as in~\cite{Cat12},
\[
\psi(x) = \begin{cases}
-\log(1-x+x^2/2),~x<0 \\
\log(1+x+x^2/2),~x\geq 0.
\end{cases}
\]
Let~$(\lambda_i)$ denote a sequence of positive numbers. Choose positive sequences~$(a_n)$ and $(b_n)$ with~$a_n \leq b_n$ for each~$n$ and define
\begin{align} \label{eq:LU}
    L_n(X_1,\cdots,X_n) &:= \text{solution of the equation}~\sum_{i = 1}^n \psi(\lambda_i (X_i - x)) = b_n, \nonumber\\
    U_n(X_1,\cdots,X_n) &:= \text{solution of the equation}~\sum_{i = 1}^n \psi(\lambda_i (X_i - x)) = a_n.
\end{align}
\cite[Theorem 9 \& 10]{WR22} shows that, for a specific choice of the sequences $(a_n)$ and $(b_n)$, \ref{eq:LU} defines a confidence sequence for~$\mu$. Furthermore, under certain conditions on the sequence~$(\lambda_i)$, for all~$n$ large enough,
\begin{equation} \label{eq:cbd_wr}
    W_n(X_1,\cdots,X_n) \leq 4 \frac{\sigma^2 \sum_{i=1}^n \lambda_i^2 + G(\alpha,\delta)}{\sum_{i=1}^n \lambda_i}
\end{equation}
for some~$n-$independent constant~$G(\alpha,\delta)$ and~$\delta \in (0,1)$. We are interested in establishing a lower bound nearly matching~\eqref{eq:cbd_wr}. We will allow any choice of the sequences $(a_n)$ and $(b_n)$ and a larger flexibility in the choice of the function~$\psi$.

\begin{theorem} \label{thm:LB}
Let~$\texttt{Var}(X_1) = \sigma^2$ and $\mathbb{E}|X_1|^{2+\vartheta} < \infty$ for some~$0 < \vartheta \leq 1$. Let the sequence~$(\lambda_i)$ be such that $\lambda_i \downarrow 0$ and $\sum_{i=1}^\infty \lambda_i^2 = \infty$. Let the influence function~$\psi$ be $1-$Lipschitz. Then, the width of the confidence sequence~$W_n(X_1,\cdots,X_n)$ is such that
\[
W_n(X_1,\cdots,X_n) \geq \frac{a \Big(\sum_{i=1}^n \lambda_i^2 \log \log \sum_{i=1}^n  \lambda_i^2\Big)^{1/2}}{\sum_{i=1}^n \lambda_i},
\]
for any~$a < 2 \sigma \sqrt{2}$ and at least for large~$n$.
\end{theorem}
The result provides the minimum width of the confidence sequence in the finite variance case. Below, we sketch the broad ideas employed in establishing the result. Notice that since~the influence function~$\psi$ is $1-$Lipschitz, we have~$\Big|\frac{d}{dx} \sum_{i=1}^n \psi(\lambda_i(X_i-x)) \Big| \leq \sum_{i=1}^n \lambda_i.$
Therefore, the confidence intervals defined in \eqref{eq:LU} satisfy
\[
W_n(X_1,\cdots,X_n) \geq  \frac{b_n - a_n}{\sum_{i=1}^n \lambda_i}.
\]
Next, if \eqref{eq:LU} defines the confidence sequence for~$\mu$ we need to show that~$b_n - a_n$ cannot be too small. Indeed, from \eqref{eq:LU} we have
\begin{equation} \label{eq:p_cs}
\mathbb{P}\Big(a_n \leq \sum_{i=1}^n \psi(\lambda_i (X_i-m)) \leq b_n~\text{for all}~n \geq 1 \Big) \geq 1 - \alpha.
\end{equation}
Consider the transformation~$Y_i = \psi(\lambda_i(X_i-m))$,
whence $Y_1,Y_2,\cdots$ is a sequence of independent random variables with all finite moments. For~$n\geq 1$, let~$\widetilde{a}_n = a_n - \sum_{i=1}^n \mathbb{E}Y_i,~~\widetilde{b}_n = b_n -\sum_{i=1}^n \mathbb{E}Y_i.$
Then~\eqref{eq:p_cs} implies that
\begin{equation} \label{eq:p_cs_t}
    \mathbb{P}\Big(\widetilde{a}_n \leq \sum_{i=1}^n (Y_i - \mathbb{E}Y_i) \leq \widetilde{b}_n~\text{for all}~n \geq 1 \Big) \geq 1 - \alpha.
\end{equation}
Clearly,~$\widetilde{b}_n - \widetilde{a}_n = b_n - a_n$, and we will show that these differences cannot be too small using~\eqref{eq:p_cs} and the general law of iterated logarithm in~\cite{Wit85}.

\vspace{-0.05in}
\section{Catoni-style Confidence Sequence for~$p \in (1,2]$} \label{sec:IV}
\vspace{-0.05in}

From Theorem~\ref{thm:LB}, it is clear that the lower bounds are sharper than the width obtained (see \eqref{eq:cbd_wr}) in~\cite{WR22}. In this section, we extend the analysis and establish the width of the confidence interval for~$p = (1,2]$, and strengthen the result in~\cite{WR22} for~$p=2$. Note that the constant~$G(\alpha,\delta)$ in \eqref{eq:cbd_wr} depends on an additional parameter~$\delta \in (0,1)$, with the width holding with probability~$1-\delta$. We both extend the approach in~\cite{WR22} to the case of
infinite variance and strengthen their result, even in the case of
finite variance, where the constant is only a function of~$\alpha$.

For an appropriate $C_p>0$ let $\psi$  be a non-decreasing
 function $\bbr\to\bbr$ such that
\begin{equation} \label{e:psi.catoni.p}
  -\log(1-x+C_p|x|^p)\leq \phi(x)\leq \log(1+x+C_p|x|^p)
\end{equation}
for all $x\in\bbr$.  One way to chose~$C_p = \Big(\frac{p-1}{p} \Big)^{p/2} \Big(\frac{2-p}{p-1} \Big)^{(2-p)/2}$, whence~$C_2 = 1/2$ as in~\cite{Cat12}. Let $(\lambda_n)$ be a sequence of positive numbers such that
\begin{equation} \label{e:lambda}
  \lim_{n\to\infty} \lambda_n=0, \ \ \ \sum_{n=1}^\infty \lambda_n^p=\infty.
\end{equation}
It follows immediately that the processes~$M_n^+ = \prod_{i=1}^n \exp\bigl\{ \phi\bigl( \lambda_i(X_i-\mu)\bigr) -
  C_pv_p\lambda_i^p\bigr\}$ and $M_n^- = \prod_{i=1}^n \exp\bigl\{ -\phi\bigl( \lambda_i(X_i-\mu)\bigr) -
  C_pv_p\lambda_i^p\bigr\}$~are non-negative supermartingales with respect to the natural
filtration of the sequence $X_1,X_2,\ldots$. Let $0<\alpha<1$ be a
confidence level. By the maximal inequality for non-negative
supermartingales, the following sequence of sets forms a
$(1-\alpha)$-confidence sequence for $\mu$:
\begin{align} \label{e:conf.seq}
I_n(\alpha)=\left\{ x\in\bbr:\, -\log\frac{2}{\alpha}-C_pv_p\sum_{i=1}^n
  \lambda_i^p \leq \sum_{i=1}^n \phi\bigl( \lambda_i(X_i-x)\bigr)
  \leq \log\frac{2}{\alpha}+C_pv_p\sum_{i=1}^n
  \lambda_i^p \right\},
\end{align}
$n=1,2,\ldots$. Note that by \eqref{e:lambda}, the sum $\sum_{i=1}^n
  \lambda_i^p $ grows slower than linearly fast with $n$. Therefore,
  at least for large $n$ the equations
  \begin{equation} \label{e:end.points}
 \sum_{i=1}^n \phi\bigl( \lambda_i(X_i-x)\bigr)
 =\pm\left(\log\frac{2}{\alpha}+C_pv_p\sum_{i=1}^n \lambda_i^p\right)
 \end{equation}
has unique real roots, which we will denote by $x_{+,n}$  and $x_{-,n}$ respectively, in which case the set
$I_n(\alpha)$ is an interval of the finite length~$\bigl| I_n(\alpha)\bigr| =x_{-,n}-x_{+,n}.$
The next result characterizes how fast these lengths grow as~$n$ increases.

\begin{theorem} \label{thm:CS_iv}
Suppose the sequence~$(\lambda_n)$ is non-random,~$0< t_n < 1$ and~$\tau_n > 0$. Suppose~$\varepsilon_n = \alpha \exp\Big\{-C_p \upsilon_p \sum_{i=1}^n \lambda_i^p (1+ t_i)^{-(p-1)} \Big\}$ for~$n = 1,2,\cdots$. Consider the condition
\begin{align}  \label{eq:cond_mr}
&C_pv_p \sum_{i=1}^n \lambda_i^p(1+t_i^{-(p-1)}) +\log 2/\alpha
+\log 2/\vep_n  \leq
\frac{\tau_n^{1/(p-1)}}{(1+\tau_n)^{p/(p-1)}}
\frac{\left(\sum_{i=1}^n \lambda_i\right)^{p/(p-1)}}{\left(C_p \sum_{i=1}^n
  \lambda_i^p (1-t_i)^{-(p-1)}\right)^{1/(p-1)}}.
  \end{align}
   The following holds for the width of the confidence sequence
 \begin{equation*}
     P\biggl( \bigl| I_n(\alpha)\bigr|\leq
4(1+\tau_n) \frac{C_pv_p\sum_{i=1}^n \lambda_i^p\bigl(1+t_i^{-(p-1)}\bigr)
   +\log 2/\alpha}{\sum_{i=1}^n
  \lambda_i}, \text{~$\forall n$ such that~\eqref{eq:cond_mr} holds} \biggr) \geq 1 - \alpha \sum_{i=1}^{\infty} \varepsilon_i.
 \end{equation*}
\end{theorem}

Note that with $\vep_n$ chosen as in the statement, the condition \eqref{eq:cond_mr} holds for all large $n$, at least if~$(t_n)$ are bounded away from 0, and $(\tau_n)$ are not too small (we
could, in fact, keep $\tau_n$ a small positive constant). The right hand side can be made small if $\alpha$ is small.

Note that~$W_n(X_1,\cdots,X_n):=\bigl| I_n(\alpha)\bigr|$, so comparing with~\eqref{eq:cbd_wr} for~$p=2$, we note the following differences: (i)~The width in~\eqref{eq:cbd_wr} depends on both~$\alpha$ and another confidence parameter~$\delta$ implying a compromise in the width, while the width in Theorem~\ref{thm:CS_iv} depends only on the confidence parameter~$\alpha$. (ii)~Abusing the notation, let~$W_{n,\lambda} = W_n(X_1,\cdots,X_n) \sum_{i=1}^n \lambda_i$. For~$p=2$, using Theorem~\ref{thm:CS_iv} ~$W_{n,\lambda} \leq 2(1+\tau_n) \sigma^2 \sum_{i=1}^n \lambda_i^2 (1+1/t_i) + \log(2/\delta)$, while that from~\eqref{eq:cbd_wr} is~$W_{n,\lambda} \leq 4 \sigma^2 \sum_{i=1}^n \lambda_i^2 + \log(2/\delta) + \log(2/\alpha)$; implying sharper bounds from Theorem~\ref{thm:CS_iv}.

Also, $\lambda_i = \Theta(1/t^{1/p})$ from~\cite{Cat12,CJLX21} implying that the Catoni-style confidence sequence enjoys~$O(\frac{\log t}{t^{\frac{p-1}{p}}})$ shrinkage rate. The dependence on~$\alpha$ is~$O(\log(1/\alpha))$ improving over~Theorem~\ref{thm:DS}, as the width increases slowly in case of Catoni-style sequence as~$\alpha \downarrow 0$.

\vspace{-0.05in}
\section{Proofs of Main Results}
\vspace{-0.05in}

In this section, we provide proofs of the main results to illustrate the challenges involved.

\textbf{Proof of Theorem~\ref{thm:DS}}: Let~$S^+_n = \sum_{i=1}^n \lambda_i (X_i-\mu)$ and $S^-_n = \sum_{i=1}^n - \lambda_i (X_i-\mu)$ denote two martingales. The confidence intervals and hence the sequence is obtained by applying Lemma~\ref{lem:DS} to each of these martingales.
Let~$a = \frac{1}{m_p b^{\frac{1}{p-1}}} \Big( \Big( \frac{2}{\alpha}\Big)^{\frac{1}{p-1}} - 1\Big)$. We have from Lemma~\ref{lem:DS},
\begin{align*}
    \mathbb{P}\Bigg( \forall~n, \sum_{i=1}^n \lambda_i (X_i - \mu)  \leq a + b \sum_{i=1}^n \lambda_i^2 \mathbb{E}[|X_i-\mu|^p|\mathcal{F}_{n-1}] \Bigg) &\geq 1 - \alpha/2, \\
    \mathbb{P}\Bigg( \forall~n, -\sum_{i=1}^n \lambda_i (X_i - \mu)  \leq a + b \sum_{i=1}^n \lambda_i^2 \mathbb{E}[|X_i-\mu|^p|\mathcal{F}_{n-1}] \Bigg) &\geq 1 - \alpha/2.
 \end{align*}
 By using the fact that~$\mathbb{E}[|X_i-\mu|^p|\mathcal{F}_{n-1}] \leq \upsilon_p$ and taking an union bound the result follows. The sequence that optimizes the width is calculated using~\cite[Eq. (24-28)]{WR20} and \cite[Appendix A]{WR22} as the minimizer of~$b \upsilon_p \lambda^{p-1} + \frac{a}{t \lambda},$
solving which we obtain~the desired sequence.

\textbf{Proof of Theorem~\ref{thm:LB}}:
Let~$s^2_n = \sum_{i=1}^n \texttt{Var}(Y_i)$ for~$n = 1,2,\cdots$. We will first show that~$\texttt{Var}(Y_i) \sim \lambda_i^2 \sigma^2$ as~$i \rightarrow \infty$. Indeed, $\mathbb{E}Y_i = \mathbb{E}\psi(\lambda_i(X - \mu)) \leq \mathbb{E}[\log (1+ \lambda_i (X - \mu)) + \lambda_i^2 (X-\mu)^2/2] \leq \mathbb{E}[\lambda_i(X-\mu) + \lambda_i^2 (X-\mu)^2/2] = \lambda_i^2 \sigma^2/2$. Also, $\mathbb{E}Y_i \geq \mathbb{E}[- \log(1 - \lambda_i(X-\mu) + \lambda_i^2(X-\mu)^2)/2] \geq \mathbb{E}[-(-\lambda_i (X-\mu) + \lambda_i^2 (X-\mu)^2/2)] = - \lambda_i^2 \sigma^2/2.$ Therefore, $|\mathbb{E}Y_i| \leq \lambda_i^2 \sigma^2/2$ for $i = 1,2,\cdots$. Next,
\begin{align*}
\mathbb{E}Y_i^2 = \mathbb{E}\psi^2(\lambda_i(X-\mu)) &\leq \mathbb{E}\Big\{ [\log(1+\lambda_i(X-\mu) + \lambda_i^2(X-\mu)^2)/2]^2 \mathds{1}(X \geq \mu) \Big\} \\ &+  \mathbb{E}\Big\{ [\log(1-\lambda_i(X-\mu) + \lambda_i^2(X-\mu)^2)/2]^2 \mathds{1}(X < \mu) \Big\}.
\end{align*}
There is~$x_0>0$ such that~$\log(1+x) \leq x^{1/2}$ for $x \geq x_0$. We have
\begin{align} \label{eq:up_bd}
    & \mathbb{E}\Big\{ [\log(1+\lambda_i(X-\mu) + \lambda_i^2(X-\mu)^2)/2]^2 \mathds{1}(X \geq \mu) \Big\} \nonumber\\ &\leq \mathbb{E}\Big\{ [\lambda_i (X-\mu) + \lambda_i^2(X - \mu)^2/2]^2 \mathds{1}(\mu \leq X \leq \mu + x_0/\lambda_i \Big\} \nonumber\\
    &+\mathbb{E}\Big\{ [\lambda_i (X-\mu) + \lambda_i^2(X - \mu)^2/2]^2 \mathds{1} (X > \mu + x_0/\lambda_i)\Big\} \nonumber\\
    &= \lambda_i^2 \mathbb{E}[(X - \mu)^2 \mathds{1} (\mu \leq X \leq \mu + x_0/\lambda_i)] + o(\lambda_i^2) \nonumber\\
    &= \lambda_i^2 \mathbb{E}[(X-\mu)^2 \mathds{1}(X \geq \mu)] + o(\lambda_i^2).
\end{align}
Similarly,
\begin{align} \label{eq:low_bd}
    &\mathbb{E}\Big\{ [\log(1-\lambda_i(X-\mu) + \lambda_i^2(X-\mu)^2)/2]^2 \mathds{1}(X < \mu) \Big\} =\lambda_i^2 \mathbb{E}[(X-\mu)^2 \mathds{1}(X < \mu)] + o(\lambda_i^2).
\end{align}
From~\eqref{eq:up_bd} and \eqref{eq:low_bd}, we have that~$\mathbb{E}Y_i^2 \leq \lambda_i^2 \sigma^2 + o(\lambda_i^2)$. On the other hand,
\begin{align*}
\mathbb{E}Y_i^2  &\geq \mathbb{E}\Big\{ [\log(1-\lambda_i(X-\mu) + \lambda_i^2(X-\mu)^2)/2]^2 \mathds{1}(\mu \leq X \leq \mu + x_0/\lambda_i) \Big\} \\ &+  \mathbb{E}\Big\{ [\log(1+\lambda_i(X-\mu) + \lambda_i^2(X-\mu)^2)/2]^2 \mathds{1}(\mu - x_0/\lambda_i \leq X < \mu) \Big\}.
\end{align*}
Similarly, $\mathbb{E}Y_i^2  \geq \lambda_i^2 \sigma^2 + o(\lambda_i^2)$. Therefore,
\[
\mathbb{E}Y_i^2  = \lambda_i^2 \sigma^2 + o(\lambda_i^2).
\]
It follows from the above arguments that
\[
s_n^2 \sim \sigma^2 \sum_{i=1}^n \lambda_i^2,~\text{and}~\theta_n:= s_n (2 \log \log s_n^2)^{1/2}~ \sim~ \sigma \Big( 2 \sum_{i=1}^n \lambda_i^2 \log \log \sum_{i=1}^n \lambda_i^2\Big)^{1/2}.
\]
We verify that the condition~$(2,\alpha)$ in~\cite{Wit85} holds for sequence~$\theta_n$ and $\vartheta$. Denoting by~$c$ a generic positive constant that may change from time to time, we have for large~$n_0$,
\[
\sum_{n=n_0}^{\infty} \theta_n^{-(2+\vartheta)} \mathbb{E}|Y_i - \mathbb{E}Y_i|^{2+\vartheta} \leq c \sum_{n=n_0}^{\infty} \Big( \sum_{i=1}^n \lambda_i^2 \log \log \sum_{i=1}^n \lambda_i^2 \Big)^{-1 - \vartheta/2} \mathbb{E}|Y_i|^{2+\vartheta}.
\]
We have using~$\log(1+x+x^2/2) \leq 2 \log(1+x/\sqrt{2})$ and $\log(1-x+x^2/2) \geq 2 \log(1-x/\sqrt{2})$,
\begin{align*}
    \mathbb{E}|Y_i|^{2+\vartheta} &\leq \mathbb{E}[[2 \log(1+\lambda_i(X-\mu)/\sqrt{2})]^{2+\vartheta} \mathds{1}(X > \mu)] \nonumber \\
    &+ \mathbb{E}[[2 \log(1-\lambda_i(X-\mu)/\sqrt{2})]^{2+\vartheta} \mathds{1}(X < \mu)] \leq c \lambda_i^{2+\vartheta}.
\end{align*}
Therefore, as ($\lambda_n$) is non-increasing,
\[
\sum_{n=n_0}^{\infty} \theta_n^{-(2+\vartheta)} \mathbb{E}|Y_i - \mathbb{E}Y_i|^{2+\vartheta} \leq c \sum_{n=n_0}^{\infty} \Big( \sum_{i=1}^n \lambda_i^2  \Big)^{-1 - \vartheta/2} \lambda_i^{2+\vartheta} \leq c \sum_{n=n_0}^{\infty} \frac{1}{n^{1+\vartheta/2}} < \infty.
\]
Hence~the condition~$(2,\alpha)$ in~\cite{Wit85} holds. It follows that
\[
\limsup_{n \rightarrow \infty} \frac{s_{n+1}}{s_n} = \limsup_{n \rightarrow \infty} \Big(\frac{\sum_{i=1}^{n+1} \lambda_i^2}{ \sum_{i=1}^n \lambda_i^2} \Big)^{1/2} \leq \limsup_{n \rightarrow \infty} \Big(\frac{n+1}{n} \Big)^{1/2} = 1.
\]
Therefore, by \cite[Theorem 2.1]{Wit85}
\[
\limsup_{n \rightarrow \infty} \theta_n^{-1} \sum_{i=1}^n (Y_i - \mathbb{E}Y_i) = 1,~\text{a.s.}
\]
That is,
\[
\limsup_{n \rightarrow \infty} \frac{\sum_{i=1}^n (Y_i - \mathbb{E}Y_i)}{\sigma \Big( 2 \sum_{i=1}^n \lambda_i^2 \log \log \sum_{i=1}^n \lambda_i^2\Big)^{1/2}} = 1,~\text{a.s.}
\]
This means that the sequence~$(\widetilde{a}_n)$ and $(\widetilde{b}_n)$ must satisfy
\[
\widetilde{b}_n - \widetilde{a}_n > a \Big( 2 \sum_{i=1}^n \lambda_i^2 \log \log \sum_{i=1}^n \lambda_i^2\Big)^{1/2}
\]
for any~$a < 2 \sigma \sqrt{2}$ and large~$n$, which implies the same for~$b_n - a_n$ and the result holds.

\textbf{Proof of Theorem~\ref{thm:CS_iv}}: For a fixed $x\in\bbr$ the following processes are also non-negative
supermartingales:
\begin{align} \label{e:Mx.plus}
&M_n^+(x) =  \prod_{i=1}^n \exp\bigl\{ \phi\bigl(
  \lambda_i(X_i-x)\bigr)\bigr\} \\
\notag &\hskip 0.7in \exp\left\{
  -(\mu-x)\sum_{i=1}^n \lambda_i -
  C_pv_p\sum_{i=1}^n \lambda_i^pt_i^{-(p-1)}-C_p|\mu-x|^p
  \sum_{i=1}^n \lambda_i^p(1-t_i)^{-(p-1)}\right\}
\end{align}
and
\begin{align} \label{e:Mx.minus}
&M_n^-(x) =  \prod_{i=1}^n \exp\bigl\{ -\phi\bigl(
  \lambda_i(X_i-x)\bigr)\bigr\} \\
\notag &\hskip 0.7in \exp\left\{
  (\mu-x)\sum_{i=1}^n \lambda_i -
  C_pv_p\sum_{i=1}^n \lambda_i^pt_i^{-(p-1)}-C_p|\mu-x|^p
  \sum_{i=1}^n \lambda_i^p(1-t_i)^{-(p-1)}\right\} .
\end{align}
Note that for $x=\mu$ and $t_i\equiv 1$, these processes
become $M_n^+$ and $M_n^-$, correspondingly. Denote~$f_n(x)=\sum_{i=1}^n \phi\bigl(
\lambda_i(X_i-x)\bigr).$
The maximal inequality for non-negative supermartingales, for every
$x\in\bbr$ and $h>0$,
\begin{align*}
&P\left(\exp\left\{ f_n(x) -(\mu-x)\sum_{i=1}^n \lambda_i -
  C_pv_p\sum_{i=1}^n \lambda_i^pt_i^{-(p-1)}-C_p|\mu-x|^p
                 \sum_{i=1}^n \lambda_i^p(1-t_i)^{-(p-1)}\right\} \geq h\right)\\
\leq&1/h,
\end{align*}
which is the same as
\begin{align} \label{e:fn.bound1}
P\biggl( &f_n(x) \geq (\mu-x)\sum_{i=1}^n \lambda_i +
  C_pv_p\sum_{i=1}^n \lambda_i^pt_i^{-(p-1)}+C_p|\mu-x|^p
                 \sum_{i=1}^n \lambda_i^p(1-t_i)^{-(p-1)}
  + \log h \biggr) \leq 1/h.
\end{align}

Choose $h=2/\vep_n$ for $0<\vep_n<1$ and denote
\begin{align} \label{e:Bplus}
B_n^+(x) = (\mu-x)\sum_{i=1}^n \lambda_i +
  C_pv_p\sum_{i=1}^n \lambda_i^pt_i^{-(p-1)}+C_p|m-x|^p
                 \sum_{i=1}^n \lambda_i^p(1-t_i)^{-(p-1)}
  + \log 2/\vep_n.
\end{align}
Then \eqref{e:fn.bound1} translates into
\begin{align} \label{e:fn.bound2}
P\Bigl( f_n(x) \geq B_n^+(x)\Bigr) \leq
  \vep_n/2.
\end{align}

Consider now the equation
\begin{equation} \label{e:eq.plus}
  B_n^+(x) = -C_pv_p \sum_{i=1}^n \lambda_i^p -\log 2/\alpha.
\end{equation}
We will establish  conditions under which this equation
has real roots. Assuming, for a moment, that such roots exist, let
$y_n$ denote the smallest such root. Using \eqref{e:eq.plus} with
$x=y_n$ tells us that on an event of probability at least
$1-\vep_n/2$,~we have $
f_n(y_n)<-C_pv_p \sum_{i=1}^n \lambda_i^p -\log 2/\alpha.
$
We conclude by the definition of $x_{-,n}$ in \eqref{e:conf.seq}, that
\begin{equation} \label{e:comp.xy}
P\Bigl( x_{-,n}<y_n \ \ \text{for all $n$ for which \eqref{e:eq.plus}
  has real roots} \Bigr) \geq
1-\sum_{n=1}^\infty \vep_n/2.
\end{equation}

We now establish conditions for the equation \eqref{e:eq.plus}
  to have real roots. The function $B_n^+$ is a strictly convex
  function of $x$, diverging to infinity at $\pm\infty$, so it has a
  unique minimum, achieved at the point
  $$
  z_+= \mu+\left( \frac{\sum_{i=1}^n \lambda_i}{pC_p \sum_{i=1}^n
      \lambda_i^p (1-t_i)^{-(p-1)}}\right)^{1/(p-1)},
  $$
  and we have
  \begin{equation} \label{e:B.min}
  B_n^+(z_n)=-\frac{p}{p-1} \frac{\left(\sum_{i=1}^n \lambda_i\right)^{p/(p-1)}}{\left(pC_p \sum_{i=1}^n
      \lambda_i^p (1-t_i)^{-(p-1)}\right)^{1/(p-1)}}
  + C_p v_p\sum_{i=1}^n
      \lambda_i^p t_i^{-(p-1)}+\log 2/\vep_n.
\end{equation}
If this minimal value satisfies
\begin{align} \label{e:cond1.sol}
&-\frac{p}{p-1} \frac{\left(\sum_{i=1}^n \lambda_i\right)^{p/(p-1)}}{\left(pC_p \sum_{i=1}^n
      \lambda_i^p (1-t_i)^{-(p-1)}\right)^{1/(p-1)}}
  + C_p v_p\sum_{i=1}^n
  \lambda_i^p t_i^{-(p-1)}+\log 2/\vep_n \\
 \notag \leq &-C_pv_p \sum_{i=1}^n \lambda_i^p -\log 2/\alpha,
  \end{align}
then  the equation \eqref{e:eq.plus} has real roots. Note that we can
rewrite the condition \eqref{e:cond1.sol} in the form
\begin{align} \label{e:cond2.sol}
C_pv_p \sum_{i=1}^n \lambda_i^p(1+t_i^{-(p-1)}) +\log 2/\alpha
+\log 2/\vep_n \leq
\frac{p}{p-1} \frac{\left(\sum_{i=1}^n \lambda_i\right)^{p/(p-1)}}{\left(pC_p \sum_{i=1}^n
  \lambda_i^p (1-t_i)^{-(p-1)}\right)^{1/(p-1)}}.
  \end{align}
We claim that this condition holds for all large $n$, at least if
$(t_n)$  are bounded away from 0, and if $\vep_n$ is not too small.
Indeed, in this case for some
constant $C$,
\begin{equation} \label{e:num}
\sum_{i=1}^n \lambda_i^p(1+t_i^{-(p-1)})
\leq C \sum_{i=1}^n \lambda_i^p,
\end{equation}
while
\begin{equation} \label{e:den}
\frac{\left(\sum_{i=1}^n \lambda_i\right)^{p/(p-1)}}{\left(  \sum_{i=1}^n
    \lambda_i^p (1-t_i)^{-(p-1)}\right)^{1/(p-1)}} \geq
\frac{\left(\sum_{i=1}^n \lambda_i\right)^{p/(p-1)}}{\left(  \sum_{i=1}^n
    \lambda_i^p\right)^{1/(p-1)}} .
\end{equation}
Since the ratio of the expressions in the right-hand sides of
\eqref{e:den} and \eqref{e:num} is
$$
\frac1C \left( \frac{\left(\sum_{i=1}^n \lambda_i\right)}{\sum_{i=1}^n
    \lambda_i^p}\right)^{1/(p-1)} \to \infty
$$
by \eqref{e:lambda}, we conclude that \eqref{e:cond2.sol} holds and,
hence, the equation \eqref{e:eq.plus}   has real roots, at
least for all large $n$, as long $\vep_n$ does not go to zero too
fast.

Notice that, if $\vep_n\leq 2$, then
$$
B_n^+(\mu) = C_pv_p\sum_{i=1}^n \lambda_i^pt_i^{-(p-1)}+\log 2/\vep_n
>0 > -C_pv_p \sum_{i=1}^n \lambda_i^p -\log 2/\alpha.
$$
Furthermore, the minimum of $B_n^+$ is achieved to the right of
$\mu$. Therefore, under the condition \eqref{e:cond2.sol},  the equation
\eqref{e:eq.plus} has one or two real roots to the right of $\mu$, and
$y_n$ is the smallest of these roots.

For $x>\mu$ the equation \eqref{e:eq.plus} becomes
\begin{align} \label{e:right.of.m}
&C_p(x-\mu)^p \sum_{i=1}^n \lambda_i^p(1-t_i)^{-(p-1)}
- (x-\mu)\sum_{i=1}^n \lambda_i\\\notag
&\hspace{1.5in}+C_pv_p\sum_{i=1}^n \lambda_i^p\bigl(1+t_i^{-(p-1)}\bigr)
  + \log 2/\vep_n+\log 2/\alpha=0.
\end{align}
We can rewrite
\eqref{e:right.of.m} in the form
\begin{equation} \label{e:right.of.m1}
  Kz^p-z+M=0
\end{equation}
for $z=x-\mu>0$ and
$$
K=\frac{C_p \sum_{i=1}^n \lambda_i^p(1-t_i)^{-(p-1)}}{\sum_{i=1}^n
  \lambda_i}
$$
and
$$
M= \frac{C_pv_p\sum_{i=1}^n \lambda_i^p\bigl(1+t_i^{-(p-1)}\bigr)
  + \log 2/\vep_n+\log 2/\alpha}{\sum_{i=1}^n
  \lambda_i}.
$$
Setting $y=K^{1/(p-1)}z>0$ and $D=K^{1/(p-1)}M$ transforms
\eqref{e:right.of.m1} into the equation
\begin{equation} \label{e:right.of.m2}
  y^p-y+D=0.
\end{equation}
Let $\tau_n>0$ and suppose that
\begin{equation} \label{e:cond.D}
  D\leq \frac{\tau_n^{1/(p-1)}}{(1+\tau_n)^{p/(p-1)}}.
  \end{equation}
Then the equation \eqref{e:right.of.m2} has a positive solution $y(D)$
satisfying
$$
y(D)\leq (1+\tau_n)D,
$$
which implies that
\begin{equation} \label{e:yn.bound}
  y_n\leq \mu+(1+\tau_n)M =
  \mu+(1+\tau_n) \frac{C_pv_p\sum_{i=1}^n \lambda_i^p\bigl(1+t_i^{-(p-1)}\bigr)
  + \log 2/\vep_n+\log 2/\alpha}{\sum_{i=1}^n
  \lambda_i}.
\end{equation}
Note that the condition \eqref{e:cond.D} can be rewritten in the form
\begin{align} \label{e:condD.sol}
&C_pv_p \sum_{i=1}^n \lambda_i^p(1+t_i^{-(p-1)}) +\log 2/\alpha
+\log 2/\vep_n  \leq
\frac{\tau_n^{1/(p-1)}}{(1+\tau_n)^{p/(p-1)}}
\frac{\left(\sum_{i=1}^n \lambda_i\right)^{p/(p-1)}}{\left(C_p \sum_{i=1}^n
  \lambda_i^p (1-t_i)^{-(p-1)}\right)^{1/(p-1)}}.
  \end{align}
Similarly to the condition \eqref{e:cond2.sol}, this condition holds
for all large $n$ as long as $\vep_n$ and $\tau_n$ do not go to zero
too fast. We conclude by \eqref{e:comp.xy} and \eqref{e:yn.bound}
\begin{align} \label{e:comp1.xy}
&P\biggl( x_{-,n}<\mu+(1+\tau_n) \frac{C_pv_p\sum_{i=1}^n \lambda_i^p\bigl(1+t_i^{-(p-1)}\bigr)
  + \log 2/\vep_n+\log 2/\alpha}{\sum_{i=1}^n
  \lambda_i} \\
  \notag &\hskip 1.5in \text{for all $n$ for which \eqref{e:condD.sol}
           holds} \biggr) \geq
1-\sum_{i=1}^\infty \vep_i/2.
\end{align}

The same argument shows that
\begin{align} \label{e:comp2.xy}
&P\biggl( x_{+,n}>\mu-(1+\tau_n) \frac{C_pv_p\sum_{i=1}^n \lambda_i^p\bigl(1+t_i^{-(p-1)}\bigr)
  + \log 2/\vep_n+\log 2/\alpha}{\sum_{i=1}^n
  \lambda_i} \\
\notag &\hskip 1.5in \text{for all $n$ for which  \eqref{e:condD.sol}
  holds} \biggr) \geq
1-\sum_{i=1}^\infty \vep_i/2.
\end{align}

We conclude by  \eqref{e:comp1.xy} and
\eqref{e:comp2.xy} that
\begin{align} \label{e:length.bound}
&P\biggl( \bigl| I_n(\alpha)\bigr|\leq
2(1+\tau_n) \frac{C_pv_p\sum_{i=1}^n \lambda_i^p\bigl(1+t_i^{-(p-1)}\bigr)
  + \log 2/\vep_n+\log 2/\alpha}{\sum_{i=1}^n
  \lambda_i} \\
\notag &\hskip 1.5in \text{for all $n$ for which  \eqref{e:condD.sol}
  holds} \biggr) \geq
1-\sum_{n=1}^\infty \vep_n.
\end{align}
With~$\varepsilon_n$ as in the statement,~\eqref{e:length.bound} is transformed into
\begin{align} \label{e:length.bound1}
&P\biggl( \bigl| I_n(\alpha)\bigr|\leq
4(1+\tau_n) \frac{C_pv_p\sum_{i=1}^n \lambda_i^p\bigl(1+t_i^{-(p-1)}\bigr)
   +\log 2/\alpha}{\sum_{i=1}^n
  \lambda_i} \\
\notag &\hskip 0.7in \text{for all $n$ for which  \eqref{e:condD.sol}
  holds} \biggr) \geq
1-\alpha \sum_{n=1}^\infty   \exp\left\{ - C_pv_p\sum_{i=1}^n
  \lambda_i^p\bigl(1+t_i^{-(p-1)}\bigr) \right\}.
\end{align}
It follows from \eqref{e:lambda} that the sum in the right hand side
is finite, and can be made small if $\alpha$ is small.

\vspace{-0.05in}
\section{Conclusion}
\vspace{-0.05in}

We provided an extension of confidence sequences for settings where
the variance of the data-generating distribution does not exist. Dealing with such heavy-tail settings required using robust estimation methods to obtain acceptable deviation bounds. We made use of the influence functions inspired by~\cite{Cat12} to obtain Catoni-style confidence sequences. We first established lower bounds on the widths of the Catoni-style confidence sequences for the finite variance case using a general law of iterated logarithm. We then provided a simple confidence sequence using a~$L_p-$version of the Dubins-Savage inequality for comparison. Finally, we derived Catoni-style confidence sequences in case of infinite variance and strengthened the existing result on finite variance in~\cite{WR22}.

\bibliographystyle{plain}
\bibliography{references}

\end{document}